\documentclass[12pt]{article}

\usepackage[latin1]{inputenc}
\usepackage{amsmath}
\usepackage{amsfonts}
\usepackage{amsthm}
\usepackage{graphicx}
\numberwithin{equation}{section}
\newtheorem{teo}{Theorem}[section]
\newtheorem{lem}[teo]{Lemma}
\newtheorem{prop}[teo]{Proposition}
\newtheorem{cor}[teo]{Corollary}
\theoremstyle{definition}

\theoremstyle{remark}
\newtheorem{rem}{Remark}

\newcommand{\ve}{\varepsilon}

\newcommand{\deb}{\rightharpoonup}

\newcommand{\impl}{\Rightarrow}

\newcommand{\pical}{\mathcal{P}}
\newcommand{\lcal}{\mathcal{L}}

\newcommand{\R}{\mathbb{R}}
\newcommand{\N}{\mathbb{N}}
\newcommand{\Z}{\mathbb{Z}}
\newcommand{\A}{\mathcal{A}}
\newcommand{\aan}{\A(\alpha,n)}
\newcommand{\clom}{\overline{\Omega}}

\newcommand{\tto}{\rightarrow}

\newcommand{\haus}{\mathcal{H}}

\newcommand\ball[2]{\overline{B(#1,#2)}}

\author{G. Buttazzo\thanks{\scriptsize\ Università di Pisa, Dip. di Matematica, Largo B. Pontecorvo, 5, 56127 Pisa, ITALY \texttt{buttazzo@dm.unipi.it}}, F. Santambrogio\thanks{\scriptsize\ Scuola Normale Superiore, Classe di Scienze, Piazza dei Cavalieri, 7, 56126 Pisa, ITALY \texttt{santambrogio@sns.it}} and N. Varchon\thanks{\scriptsize\ College Condorcet de Bresles, 60510 Bresles, FRANCE \texttt{nicolasvarchon@netscape.net}}}
\title{Asymptotics of an optimal compliance-location problem}
\date{}
\begin{document}
\maketitle
{\bf Abstract:} We consider the problem of placing a Dirichlet region made by $n$ small balls of given radius in a given domain subject to a force $f$ in order to minimize the compliance of the configuration. Then we let $n$ tend to infinity and look for the $\Gamma-$limit of suitably scaled functionals, in order to get informations on the asymptotical distribution of the centres of the balls. This problem is both linked to optimal location and shape optimization problems.

\bigskip


{\bf Math Subject Classification numbers:} 49J45, 49Q10, 74P05.

\section{Introduction}\label{sec1}

The study of asymptotical problems in the optimal location of an increasing amount of resources has been developed intensively, even in recent times, mainly by using an approach based on $\Gamma-$convergence. In \cite{BouJimRaj} the so-called location problem (choosing a set $\Sigma$ composed by $n$ points in a domain $\Omega$ in order to minimize the average distance of the points of $\Omega$ from $\Sigma$) is investigated as $n\tto +\infty$, finding a $\Gamma-$limit of a suitable sequence of functionals on the space of probability measures on $\clom$. In \cite{MosTil} the same analysis is performed for the so-called irrigation problem, where points are replaced by connected dimensional sets of finite length, and the constraint $\sharp\Sigma\leq n$ by $\haus^1(\Sigma)\leq l$. Both the problems are linked to the Monge-Kantorovich optimal transport theory. However, these asymptotical problems are not completely understood since explicit minimizing sequences are not in general known, apart some simple cases, usually in dimension two. For instance, for the location problem it is known that placing the points on a regular triangular grid, so that each one is in the middle of a cell shaped like a regular hexagon, gives an asymptotically minimizing sequence (see \cite{toth} or \cite{hexag} for stronger results).

On the other hand, many researches have been carried out on shape-optimization problems involving PDEs, i.e. optimizing the shape of a domain where to solve a PDE (usually of elliptic type with prescribed boundary conditions), in order to minimize the value of an objective functional depending on the solution of the PDE. There is a wide literature on shape optimization problems, both from a theoretical and numerical point of view. We refer for instance to the books \cite{allaire}, \cite{bensig}, \cite{bucbut}, \cite{HP}, \cite{sokzol}, where the reader can find various approaches and a lot of examples and details. It is well known (see for instance \cite{bdm1} and \cite{bdm2}) that for general cost functionals the existence of an optimal domain may fail and a relaxation procedure, involving the use of capacitary measures, is needed. It is also known that for cost functionals fulfilling particular additional
monotonicity assumptions, a simple volume constraint is sufficient to imply the existence of a minimizer (see for instance \cite{bdm3}). Moreover, these problems have shown to have many industrial applications in engineering and mechanics (bridges, light structures supporting loads\dots) and this is the reason for their widespread study and especially for numerical computations. One of the simplest shape optimization problem, which is also one of the most important in applications, is compliance minimization. It consists in finding a domain $\Omega$ (usually under a volume contraint) which minimizes the compliance value $\int_{\Omega}fu_{\Omega}\,d\lcal^d$ where $u_{\Omega}$ is the solution of the elliptic equation $-\Delta u=f$ with Dirichlet boundary conditions on $\partial\Omega$. Such a problem satisfies the monotonicity assumptions needed to have the existence of minimizers.

What we consider in this paper (see Section \ref{sec2}) is a compliance minimization problem where the unknown domain where to solve the PDE with Dirichlet boundary conditions is the complement of a finite union of balls whose number and radius are assigned. In fact, for given $n\in\N$ and $r_n>0$, we look at the problem of choosing $n$ balls $(\ball{x_i}{r_n})_{i=1,\dots,n}$ in order to minimize the compliance of $\Omega\setminus\bigcup_{i=1}^n\ball{x_i}{r_n}$ ($\Omega\subset\R^d$ and $f\in L^2(\Omega)$ are fixed). Obviously, the problem is meaningful only under a condition like $r_n\leq cn^{-1/d}$, otherwise the total volume of the balls may be sufficient to cover all $\Omega$, thus obtaining a vanishing compliance. 

The aim of the present paper is to look at the optimal location of those balls when $n\tto+\infty$. Using small balls is the most natural way to approximate a point-location problem, which would not make any sense since points have zero capacity in dimension two or more, and so they do not affect the solution of elliptic equations. In this paper we just consider the case $r_n=\alpha n^{-1/d}$, for a fixed small parameter $\alpha$, which corresponds to a volume constraint. Moreover, in Section \ref{sec5} we deal with the one dimensional case, getting it as a particular case of the $d-$dimensional, but highlighting the case of points instead of balls as well (which corresponds to taking $\alpha=0$).

The domain $\Omega$ has to be considered as an elastic membrane and $f$ stands for the forces acting on it (for instance the different loads it has to carry). Our goal is to reinforce or support the membrane at some points (by mechanical devices, by some kind of glue\dots), choosing where to locate the support points, so that it bends as less as possible, minimizing the work of the forces, given by $\int fu\,dx$. The important issue here is that we let the number of support points increase, reducing correspondingly the effect at each point, and we look at the asymptotic density of support points. 

As it happens for the location problem, placing $n$ balls with given radii in $\Omega$ reduces to a finite dimensional variational problem as a consequence of the severe geometric constraints imposed to the admissible domains. Our asymptotic analysis is based on some techniques developed in \cite{BouJimRaj} and mainly in \cite{MosTil} and so this paper, even if dealing with a shape optimization subject, stands apart from other shape optimization papers both for the setting of the problem and for its development.%Moreover, despite we minimize compliance, which is typical in shape optimization, this is not really shape optimization and in fact this paper recalls much more \cite{BouJimRaj} and \cite{MosTil}, rather than shape optimization papers. Also the techniques we use are inspirated by \cite{MosTil}, where very powerful ideas are developed, allowing more general results than those in \cite{BouJimRaj}.

Since we solve elliptic PDEs in open domains which are obtained by removing small holes from a fixed one, this subject is in connection with the problem of homogenization in perforated domain (see for instance \cite{CioMur}). Yet, there are some important differences. First, we have an optimization problem on the holes instead of taking them as given. Consequently, no periodic assumptions on the holes is supposed, even if periodic structures are used many times in the proof of the $\Gamma-$limit result. Finally, homogenization on perforated domains is interesting when the size of the holes decreases quicker than $n^{-1/d}$, otherwise, under Dirichlet boundary conditions, the solutions trivially tend to $0$. This is the reason for introducing a scaling factor which enlarge the values of the functionals and of the solutions when $n$ increases.

An interesting question we do not consider in the present paper is how to deal with the case $r_n=\alpha_n n^{-1/d}$ with $\alpha_n\tto 0$, since this could be considered as a better approximation of the case of points. Moreover, there are further reasons to study in the future this new case. First, the case $\alpha_n\tto 0$ is more linked with the theory of perforated domains developed in \cite{CioMur}, where the radius of the ball tends to zero faster than $n^{-1/d}$. Then, this case seems to require better mathematical techniques involving the behaviour of the Dirichlet energy when radii tend to $0$, since it is not simply possible to build a recovery sequence by homogenizing a given configuration (since this would let the total volume of the balls unchanged, while here it has to tend to zero). Finally, the result we get here (Theorem \ref{teoprinc}) involves a function $\theta$ we are not able to compute explicitly, while the case $\alpha_n\tto 0$ seems to require just its asymptotic expansion near $0$.

\section{Locating balls to optimize compliance}\label{sec2}

For any open set $\Omega\subset\R^d$, $\alpha>0$ and $n\in\N$ we define:
$$\aan (\Omega)=\left\{\Sigma\subset\clom\left|\Sigma=\clom\cap\bigcup_{i=1}^n\ball{x_i}{r}\text{ for }x_i\in\Omega_r,\,r=\alpha n^{-1/d}\right.\right\},$$
where $\Omega_r$ stands for the $r-$neighbourhood of $\Omega$. Given $\Omega\subset\R^d$ and $f\in L^2(\Omega)$, for any compact set $\Sigma\subset\clom$ with positive measure let us define the function $u_{f,\Sigma,\Omega}$ as the solution in the weak sense of the problem
$$\begin{cases}-\Delta u=f&\text{ in }\Omega\setminus\Sigma\\
							u=0&\text{ in }\Sigma\cup\partial\Omega,\end{cases}$$
which means precisely $u\in H^1_0(\Omega\setminus\Sigma)$ and
\begin{equation}\label{weakpde}
\int_{\Omega}\nabla u\cdot\nabla\phi\,d\lcal^d=\int_{\Omega}f\phi\,d\lcal^d\text{ for any }\phi\in H^1_0(\Omega\setminus\Sigma).
\end{equation}
Notice that $f\geq 0$ implies $u_{f,\Sigma,\Omega}\geq 0$, by the maximum principle. For $f\geq 0$, we define the compliance functional over subsets of a given domain $\Omega$ as 
$$F(\Sigma,f,\Omega)=\int_{\Omega}fu_{f,\Sigma,\Omega}\,d\lcal^d=\int_{\Omega}|\nabla u_{f,\Sigma,\Omega}|^2\,d\lcal^d.$$
\begin{rem}
The requirement $f\geq 0$ seems to be mostly a technical assumption to simplify the proofs (it allows us to deduce some pointwise inequality by maximum principle on the solutions of Dirichlet problems): physically speaking it means considering only forces which have the same direction on the whole $\Omega$ (for instance usual gravity).
\end{rem}
What we want to do now is considering such a compliance functional on the set $\aan(\Omega)$ (notice that imposing $\Sigma\in\aan(\Omega)$ implies a volume constraint given by $|\Sigma|\leq w_d\alpha^d$ and a geometrical constraint, i.e. compelling $\Sigma$ to be composed by an assigned number of identical balls). This is our $n-$th compliance minimization problem. The following existence result holds and can be proven by standard methods, due to the very severe geometry the elements of $\aan(\Omega)$ are constrained.

\begin{teo}
For any $n\in\N$, if $\Omega$ is any bounded open subset of $\R^d$ and $f\geq 0$ belongs to $L^2(\Omega)$, the problem
\begin{equation}\label{nthprob}
\min\left\{F(\Sigma,f,\Omega)\,\left|\Sigma	\in\aan(\Omega)\right.\right\}
\end{equation}
admits a solution.
\end{teo}

Then we would like to let $n$ tend to infinity and look at the asymptotics of the problem, mainly at the distribution of the centres of the balls. Let us associate to every $\Sigma\in\aan(\Omega)$ a probability measure on $\clom$, given by $\mu_{\Sigma}=n^{-1}\sum_{i=1}^n\delta_{p(x_i)}$,
where $p:\R^d\tto \clom$ is a fixed projection of the whole space to $\clom$. The role of the projection $p$ is simply to handle the case where the centre of the ball $\ball{x_i}{r}$ lies outside $\clom$. Such a measure is an atomic measure uniformly distributed on the centres (or on their projections). Then we define a functional $F_n:\pical(\clom)\tto [0,+\infty]$ by
$$F_n(\mu)=\begin{cases}n^{2/d}F(\Sigma,f,\Omega)&\text{ if }\mu=\mu_{\Sigma},\,\Sigma\in\aan(\Omega);\\
                        +\infty           &\text{otherwise.}\end{cases}$$
The coefficient $n^{2/d}$ is a factor which is needed in order to avoid the functionals to degenerate to the trivial limit functional which vanishes everywhere. Anyway, such a coefficient does not affect the choice of the minimizers. 

We will give a $\Gamma-$convergence result for the sequence $(F_n)_n$, when endowing the space $\pical(\clom)$ with the weak* topology of probability measures. To introduce the limit functional $F$ we need to define the quantity:
\begin{equation}\label{defitheta}
\theta(\alpha):=\inf\left\{\liminf_n n^{2/d}F(\Sigma_n,1,I^d)\left|\Sigma_n\in\aan(I^d)\right.\right\},
\end{equation}
where $I^d=(0,1)^d$ is the unit cube in $\R^d$.
This quantity will play the role of the constant appearing both in \cite{MosTil} (as $\theta_{n,p}$) and in \cite{BouJimRaj} (as $C_d$). However, its dependence on $\alpha$ will be essential and, unfortunately, in general not explicit. 
%For our convenience let us also define 
%$$\theta_L(\alpha):=\inf\left\{\liminf_n n^{2/d}F(\Sigma_n,1,Q)\left| \Sigma_n\in\A(\alpha_n,n)(Q),\,\alpha_n\tto\alpha\right.\right\}.$$
It is easy to see that $\theta$ is a decreasing function on $\R^+$, which vanishes after some point. In fact if $\alpha\geq\sqrt{d}/2$, it is possible to use $n$ balls of radius $\alpha n^{-1/d}$ to build a set $\Sigma\in\aan(I^d)$ covering the whole cube $I^d$, thus getting a vanishing solution $u_{f,\Sigma,I^d}=0$ and $F(\Sigma,f, I^d)=0$. Let us call $t_1$ the first vanishing point, i.e.
$$t_1:=\inf\left\{t\in\R\left|\theta(t)=0\right.\right\}\leq \frac{\sqrt{d}}{2}.$$
We denote by $\theta^-$ and $\theta^+$ the lower and upper semicontinuous envelopes of $\theta$, respectively. They are given by 
\begin{eqnarray*}\theta^-(\alpha)&=&\sup\left\{\theta(\beta)\left|\beta>\alpha\right.\right\}\\
                 \theta^+(\alpha)&=&\inf\left\{\theta(\beta)\left|\beta<\alpha\right.\right\}.\end{eqnarray*}
It is easy to check that the following formula holds:
\begin{equation}\label{thetaL}
\theta^-(\alpha)=\inf\left\{\liminf_n n^{2/d}F(\Sigma_n,1,I^d)\left| \Sigma_n\in\A(\alpha_n,n)(I^d),\,\alpha_n\tto\alpha\right.\right\}.
\end{equation}

We may now define the candidate limit functional $F$ by setting, for $\mu\in\pical(\clom)$
\begin{equation}\label{funzlim}
F(\mu)=\int_{\Omega}\frac{f^2}{\mu_a^{2/d}}\theta^-(\alpha\mu_a^{1/d})\,d\lcal^d,
\end{equation}
where $\mu_a$ denotes the density of the absolutely continuous part of $\mu$ with respect to the Lebesgue measure. It is evident from \eqref{funzlim} that the whole behaviour of the function $\theta$ affects the minimization problem for $F$.
It is more convenient to introduce the function $g_{\alpha}$ defined by $g_{\alpha}(x)=x^{-2/d}\theta(\alpha x^{1/d})$, which is a decreasing function whose semicontinuous envelops $g_{\alpha}^+$ and $g_{\alpha}^-$ are obtained by $g_{\alpha}^+(x)=x^{-2/d}\theta^+(\alpha x^{1/d})$ and $g_{\alpha}^-(x)=x^{-2/d}\theta^-(\alpha x^{1/d})$. So we have $$F(\mu)=\int_{\Omega}f^2 g_{\alpha}^-(\mu_a)\,d\lcal^d.$$

Here is a sketch of the behaviour of $g_{\alpha}$, according to what already highlighted and to what proven in Section \ref{sec4} on the function $\theta$. In particular $g_{\alpha}$ is a convex function.

\begin{figure}[hbtp]
\begin{center}
\includegraphics[height=4cm]{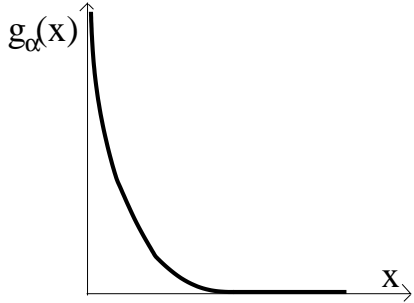}
\end{center}
\caption{Qualitative behaviour of the function $g_{\alpha}$}
\label{g1}
\end{figure}

The result we will prove is the following:
\begin{teo}\label{teoprinc}
Given any bounded open set $\Omega\subset\R^d$, a non-negative function $f\in L^2(\Omega)$ and $\alpha>0$, the sequence of functional $(F_n)_n$ previously defined $\Gamma-$converges towards $F$ with respect to the weak* topology on $\pical(\clom)$.
\end{teo}

The consequences of such a $\Gamma-$convergence result, by means of the general theory (see \cite{introgammaconve}), are the following:
\begin{itemize}
\item for any sequence $(\Sigma_n)_n$ of optimal sets for the minimization problem \eqref{nthprob} it holds, up to subsequences, $\mu_{\Sigma_n}\deb\mu$ where $\mu$ is a minimizer of $F$;
\item should $F$ have a unique minimizer, we would get full convergence of the whole sequence $\mu_{\Sigma_n}$ to the unique minimizer $\mu$.
\item the sequence of the values $\inf\left\{F(\Sigma,f,\Omega)\left|\Sigma\in\aan(\Omega)\right.\right\}$ is asymptotical to $n^{-2/d}\inf \left\{F(\mu)\left|\mu\in\pical(\clom)\right.\right\}$.
\end{itemize}

It turns out that it is very important to investigate about the minimizers of the limit functional $F$, which is much related to the behaviour of $g$. 
From Section \ref{sec4} we know that $g_{\alpha}=g_{\alpha}^-$ and that $g_{\alpha}$ is a convex and strictly decreasing function (up to the point $t_1$ where $\theta=0$). If we define $t_{\alpha}=t_1/\alpha$, it holds $g_{\alpha}(x)=0$ for any $x\geq t_{\alpha}$. We will restrict our analysis to the case where $\alpha$ is sufficiently small, so that it is not possible to cover the whole $\Omega$ by $n$ balls of radius $\alpha n^{-1/d}$.
We summarize in the following statement what we know on the minimizers of $F$.
\begin{teo}\label{opticond}
If $\alpha$ is such that $\alpha<|\Omega|t_1$, any minimizer $\mu$ for $F$ is an absolutely continuous probability measure with density $\mu_a$ which satisfies
$$\mu_a(x)\in \left(-\partial g_{\alpha}\right)^{-1}\left(\frac{c}{f^2(x)}\right),$$
for a suitable constant $c>0$. In particular it holds $\mu_a\leq t_{\alpha}$.
\end{teo}
\begin{proof}
Being $g_{\alpha}$ strictly decreasing up to $t_{\alpha}$ it is straightforward that optimality implies absolute continuity: otherwise, just remove the singular part from $\mu$ and use the same mass on the absolutely continuous part, enlarging its density and strictly decreasing the value of $F$. Moreover, by using Lagrange multipliers or performing simple variations to $\mu$, it is easy to get the existence of a constant $c>0$ such that
\begin{equation}\label{cinsubd}
f^2(x)\left(-\partial g_{\alpha}(\mu_a(x))\right)\ni c,
\end{equation}
where $\partial g$ is the, possibly multivalued, subdifferential of $g$. Then we get
$$\mu_a(x)\in \left(-\partial g_{\alpha}\right)^{-1}\left(\frac{c}{f^2(x)}\right),$$
where, for a multifunction $G$, we use the notation $G^{-1}(t)=\{z\in [0,+\infty]:t\in G(z)\}$.
We sketch a possible behaviour of $-\partial g_{\alpha}$ and $(-\partial g_{\alpha})^{-1}$ in Figure \ref{g2}.
\begin{center}
\begin{figure}[h]
\includegraphics[height=5cm]{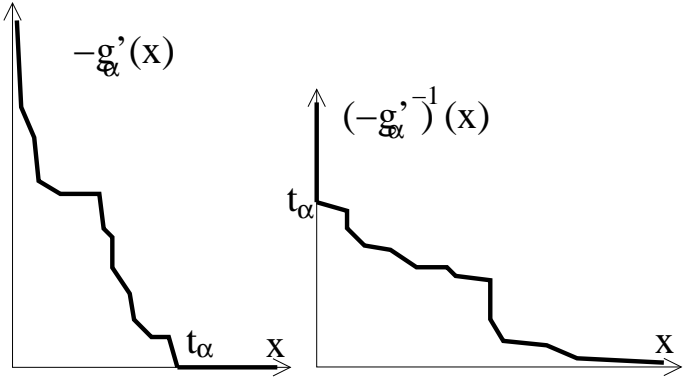}
\caption{The subdifferential $g'_{\alpha}$ of $g_{\alpha}$ and its inverse}
\label{g2}
\end{figure}
\end{center}
To get $\mu_a\leq t_{\alpha}$ (which is suggested by Figure \ref{g2} as well) it is sufficient to notice that otherwise in \eqref{cinsubd} one should have $c=0$. This would imply $\mu_a\geq t_{\alpha}$ a.e. and this is not possible by the assumption on $\alpha$, since $\mu$ has to be a probability measure.
\end{proof}

A consequence of $\mu_a\leq t_{\alpha}$ is the fact that strict convexity of $g_{\alpha}$ on $(0,t_{\alpha})$ is sufficient to ensure uniqueness of the minimizers of $F$. Section \ref{sec4} will also provide some qualitative property of $\theta$, in order to guess the behaviour of $g_{\alpha}$.

\section{The $\Gamma-$convergence result}\label{sec3}

We will prove Theorem \ref{teoprinc} in several steps, the most important two corresponding to the $\Gamma-\liminf$ and $\Gamma-\limsup$ inequalities.
\begin{prop}\label{liminf}
Under the same hypotheses of Theorem \ref{teoprinc}, denoting by $F^-$ the functional $\Gamma-\liminf_n F_n$, it holds
$F^-(\mu)\geq F(\mu)$
for any $\mu\in\pical(\clom)$. This means that, for any sequence $(\Sigma_n)_n$ such that $\mu_{\Sigma_n}$ weakly-$\star$ converges to $\mu$ and $\Sigma_n\in\aan(\Omega)$, it holds $\liminf_n n^{2/d}\int f u_n d\lcal^d\geq F(\mu)$, where $u_n$ stands for $u_{f,\Sigma_n,\Omega}$.
\end{prop}
\begin{proof}
First of all, let us fix $\ve>0$ and, in analogy to what performed in \cite{MosTil}, define the set $G_{\ve,n}$ as follows: let us assume $\Omega\subset [-a,a]^d$ for $a\in\N$, then we set
$$G_{\ve,n}=\bigcup_{y\in k^{-1}\Z^d\cap[-a,a]^d}\ball{y}{r},\,r=\alpha n^{-1/d},\,k=\left\lfloor(\ve n)^{1/d} \right\rfloor.$$
Now we define $\Sigma'_n=\Sigma_n\cup G_{\ve,n}$ and we set $u'_n=u_{f,\Sigma'_n,\Omega}$. Since $u_n\geq u'_n$, it is sufficient to estimate from below the integrals $n^{2/d}\int fu'_n d\lcal^d$. The utility of the new sequence $(u'_n)_n$ lies in the fact that $(n^{2/d}u'_n)_n$ is $L^2-$bounded. In fact we have $0\leq u'_n\leq u_{f,G_{\ve,n},\Omega}$ and, by Lemma \ref{poincare}, it holds 
$$||u_{f,G_{\ve,n},\Omega}||_{L^2(\Omega)}\leq C(\alpha,\ve,f)n^{-2/d}.$$
This implies that $(n^{2/d}u'_n)_n$ is bounded in $L^2$ and so, up to a subsequence, we have $n^{2/d}u'_n\deb w$. Since in this case we have $\liminf_n n^{2/d}\int fu_n d\lcal^d\geq \int fw d\lcal^d$, it is sufficient to estimate $w$ from below.
To do this, we first estimate the average of $w$ on a cube $Q$ centred at  point $x\in\Omega$. It holds
$$\int_Q w\,d\lcal^d=\lim_n n^{2/d}\int_Q u'_n\,d\lcal^d.$$
We use 
$$u'_n\geq u_{f,\Sigma'_n,Q}=u_{f(x),\Sigma'_n,Q}+u_{f-f(x),\Sigma'_n,Q}\geq u_{f(x),\Sigma'_n,Q}-u_{|f-f(x)|,\Sigma'_n,Q}\,\text{ in }Q,$$
where the first inequality comes from the fact that we add Dirichlet boundary conditions on $Q$ and the last by maximum principle. Notice that this is the key point in the proof where we strongly use $f\geq 0$ (in fact the other pointwise inequalities could be replaced by global integral estimates valid for general $f$, but here this is not possible, since we have to estimate integrals performed on $Q$ instead of on $\Omega$).

We will estimate separately the two terms. Let us start from the easiest, i.e. the latter. It holds $u_{|f-f(x)|,\Sigma'_n,Q}\leq u_{|f-f(x)|, G_{\ve,n},Q}$. By applying Lemma \ref{poincare} to the domain $Q$ we get
$$||u_{|f-f(x)|,\Sigma'_n,Q}||_{L^2(Q}\leq n^{-2/d}C(\alpha,\ve)||f-f(x)||_{L^2(Q)}.$$
Then we estimate, by Holder inequality,
\begin{eqnarray}\label{perlebesgue}
n^{2/d}\int_Q u_{|f-f(x)|,\Sigma'_n,Q}\,d\lcal^d&\leq& n^{2/d}|Q|^{1/2}||u_{|f-f(x)|,\Sigma'_n,Q}||_{L^2(Q}\notag\\
&\leq& C(\alpha,\ve)|Q|^{1/2}||f-f(x)||_{L^2(Q)}.
\end{eqnarray}
We now evaluate the other term. First we define the number 
$$k(n,Q)=\sharp\left(\left\{i:\, \ball{x_i}{\alpha n^{-1/d}}\cap Q\neq\emptyset\right\}\cup\left\{j:\, \ball{y_j}{\alpha n^{-1/d}}\cap Q\neq\emptyset\right\}\right).$$
Now notice that $u_{f(x),\Sigma'_n,Q}=f(x)u_{1,\Sigma'_n,Q}$. Let us denote, for simplicity, the functions $u_{1,\Sigma'_n,Q}$ by $v_n$. By a change of variables, if $\lambda$ is the side of the cube $Q$ and we define  $v_{n,\lambda}(x)=\lambda^{-2}v_n(\lambda x)$, it holds $v_{n,\lambda}=u_{1,\lambda^{-1}\Sigma'_n,I^d}$. We notice that 
$$\lambda^{-1}\Sigma'_n\in\A\left(\frac{\alpha}{\lambda}\left(\frac{k(n,Q)}{n}\right)^{1/d},k(n,Q)\right)(I^d).$$
Moreover, it holds $k(n,Q)\tto+\infty$, since 
\begin{equation}\label{ken}
k(n,Q)\geq\sharp\left\{j:\, \ball{y_j}{\alpha n^{-1/d}}\cap Q\neq\emptyset\right\}\approx \ve n |Q|.
\end{equation} 
We may also estimate the ratio between $k(n,Q)$ and $n$, by using \eqref{ken} and the fact that $\sharp\left\{i:\, \ball{x_i}{\alpha n^{-1/d}}\cap Q\neq\emptyset\right\}=\mu_n(Q_{\delta_n})$, where $\delta_n=\alpha n^{-1/d}\tto 0$ and $Q_{\delta}$ denotes the $\delta-$neighbourhood of $Q$. From $\mu_n\deb\mu$ we have $\limsup_n\mu_n(Q_{\delta_n})\leq\mu(\overline{Q})$, so that
\begin{equation}\label{kendef}
\limsup_n\frac{k(n,Q)}{n}\leq\mu(\overline{Q})+\ve|Q|.
\end{equation}
Now, by using the equality in \eqref{thetaL}, it is not difficult to derive that
\begin{eqnarray}\label{usothetaL}
\liminf_n\, k(n,Q)^{2/d}\int_{Q}v_n\,d\lcal^d&=&\liminf_n\, k(n,Q)^{2/d}\lambda^{d+2}\int_{I^d}v_{n,\lambda}\,d\lcal^d\notag\\
&\geq&\lambda^{d+2}\theta^-\left(\frac{\alpha}{\lambda}\left(\mu(\overline{Q})+\ve|Q|\right)^{1/d}\right).
\end{eqnarray}
So we get 
\begin{eqnarray*}
\liminf_n \,n^{2/d}\int_{Q}v_n\,d\lcal^d&=&\lim_n \left(\frac{n}{k(n,Q)}^{2/d}\right) \liminf_n\, k(n,Q)^{2/d}\int_{Q}v_n\,d\lcal^d\\
&\geq& \lambda^{d+2}\theta^-\left(\frac{\alpha}{\lambda}\left(\mu(\overline{Q})+\ve|Q|\right)^{1/d}\right)\left(\frac{1}{\mu(\overline{Q})+\ve|Q|}\right)^{2/d}\\
&=&\theta^-\left(\alpha\left(\frac{\mu(\overline{Q})}{|Q|}+\ve\right)^{1/d}\right)|Q|\left(\frac{|Q|}{\mu(\overline{Q})+\ve|Q|}\right)^{2/d}.
\end{eqnarray*}
This implies, recalling also \eqref{perlebesgue},
\begin{eqnarray*}
|Q|^{-1}\int_Q w\,d\lcal^d&\geq&-C(\alpha,\ve)|Q|^{-1/2}||f-f(x)||_{L^2(Q)}\\
&+&f(x)\theta^-\left(\alpha\left(\frac{\mu(\overline{Q})}{|Q|}+\ve\right)^{1/d}\right)\left(\frac{|Q|}{\mu(\overline{Q})+\ve|Q|}\right)^{2/d}.
\end{eqnarray*}
Now we let $Q$ shrink towards $x$, thus obtaining, using lower semicontinuity, for a.e. $x\in\Omega$
$$w(x)\geq f(x)\,\theta^-\left(\alpha\left(\mu_a(x)+\ve\right)^{1/d}\right)\left(\frac{1}{\mu_a(x)+\ve}\right)^{2/d},$$
where we have used the fact that a.e. point $x$ is a Lebesgue point for $f\in L^2$ to get the first term vanish, and that for any measure $\mu$ and a.e. $x\in\Omega$, it holds $|Q|^{-1}\mu(\overline{Q})\tto\mu_a(x)$ when the cube $Q$ shrinks at its centre $x$.
So we get 
\begin{eqnarray*}
\liminf_n \,n^{2/d}\int_{\Omega} fu_n d\lcal^d&\geq& \int_{\Omega} fw d\lcal^d\\
&\geq& \int_{\Omega}\frac{f^2}{(\mu_a(x)+\ve)^{2/d}}\theta^-\left(\alpha\left(\mu_a(x)+\ve\right)^{1/d}\right) d\lcal^d,
\end{eqnarray*}
and our original aim is achieved when we let $\ve\tto 0$, i.e., still using that $\theta^-$ is l.s.c.,
$$\liminf_n n^{2/d}\int_{\Omega} fu_n d\lcal^d\geq \int_{\Omega}\frac{f^2}{\mu_a(x)^{2/d}}\theta^-\left(\alpha\mu_a(x)^{1/d}\right) d\lcal^d.\qedhere$$

\end{proof}

\begin{lem}\label{poincare}
The following facts hold.
\begin{enumerate}
\item For any $0<\ve_0<1$ there exists a constant $C=C(\ve_0)$ such that
$$v\in H^1(I^d),\,|\left\{v=0\right\}|\geq \ve_0|I^d|\impl \int_{I^d} v^2d\lcal^d\leq C\int_{I^d}|\nabla v|^2d\lcal^d.$$
\item If we replace $I^d$ by a cube $Q$ whose side is $\lambda$ the same is true with the constant $\lambda^2 C$ instead of $C$.
\item As a consequence, for any $\ve>0$, any $n\in\N$, any domain $\Omega$ and $f\in L^2(\Omega)$ with $F\geq 0$, the function $u_{f,G_{\ve,n},\Omega}$, where $G_{\ve,n}$ is defined as in Theorem \ref{liminf}, satisfies $||u_{f,G_{\ve,n},\Omega}||_{L^2(\Omega)}\leq n^{-2/d}C(\alpha,\ve)||f||_{L^2(\Omega)}$.
\end{enumerate}
\end{lem}

\begin{proof}
The first assertion comes from a well-known variant of Poincar\'{e} inequality and can be proven by contradiction. The second one is just a scaling of the first. To prove the last, let us consider a family of cubes $Q_j$, each centred at a point $y_j\in\ k^{-1}\Z^d\cap[-a,a]^d$, whose side is $k^{-1}$. Let us extend the function $u_{f,G_{\ve,n},\Omega}$ to the set $\Omega'=\bigcup_j Q_j$, where the union is over the cubes touching $\Omega$, by defining a function $v$ which is identical to $u_{f,G_{\ve,n},\Omega}$ on $\Omega$ and $0$ outside (we recall that we have Dirichlet boundary conditions on $\Omega$ so that such an extension belongs to $H^1(\Omega')$). Notice that $v$ vanishes in a whole ball of radius $\alpha n^{-1/d}$ in each of these cubes, so that the ratio between the volume of such a ball and the volume of the cube depends only on $\alpha$ and $\ve$ and not on $n$. By applying the second part of the statement of this lemma, we get
$$\int_{Q_j}v^2\,d\lcal^d\leq C(\alpha,\ve)n^{-2/d}\int_{Q_j}|\nabla v|^2\,d\lcal^d,$$
and, by summing over $j$, we get
$$\int_{\Omega'}v^2\,d\lcal^d\leq C(\alpha,\ve)n^{-2/d}\int_{\Omega'}|\nabla v|^2\,d\lcal^d.$$
Since $v$ vanishes outside $\Omega$ we may write
\begin{eqnarray*}
\int_{\Omega}(u_{f,G_{\ve,n},\Omega})^2\,d\lcal^d&\leq& C(\alpha,\ve)n^{-2/d}\int_{\Omega}|\nabla u_{f,G_{\ve,n},\Omega} |^2\,d\lcal^d\\
&=&\int_{\Omega}f u_{f,G_{\ve,n},\Omega} \,d\lcal^d\leq C(\alpha,\ve)n^{-2/d}||u_{f,G_{\ve,n},\Omega}||_{L^2(\Omega)}||f||_{L^2(\Omega)}.
\end{eqnarray*}
The thesis easily follows by dividing by $||u_{f,G_{\ve,n},\Omega}||_{L^2(\Omega)}$.

\end{proof}

To get also the opposite inequality, i.e. the $\Gamma-\limsup$ inequality, we need this crucial lemma.

\begin{lem}\label{convergenza a cf}
Given $\Sigma_0\in\A(\alpha_0,n_0)(I^d)$, a domain $\Omega\subset\R^d$ and $f\in L^2(\Omega)$, we consider the sequence of sets
$$\Sigma^k=\bigcup_{y\in\ k^{-1}Z^d}\left(y+k^{-1}\Sigma_0\right)\cap\clom.$$
It happens $\Sigma^k\in\A(\alpha_0,n(k,\Omega))(\Omega)$, where $n(k,\Omega)\approx |\Omega|k^dn_0$.
Then we consider the sequence $(u_k)_k$, given by
$$u_k=k^2 u_{f,\Sigma^k,\Omega}.$$
If we assume $\partial I^d\subset\Sigma_0$, then it holds $u_k\deb c(\Sigma_0)f$ where the weak convergence is in the $L^2$ sense and $c(\Sigma_0)$ is a constant given by $\int_{I^d}u_{1,\Sigma_0,I^d}\,d\lcal^d$.
\end{lem}
\begin{proof}
First, we notice that the sequence $(u_k)_k$ is bounded in $L^2(\Omega)$. This may be proven in a way very similar to that of Lemma \ref{poincare}, but here it is even simpler. In fact we may extend as before $u_k$ to the union of the cubes of the kind $y+k^{-1}I^d$ which intersect $\Omega$, with $y\in k^{-1}\Z^d$, by giving it the value $0$ outside $\Omega$. Then we apply standard Poincaré inequality to each cube $Q_y$ ($u_k$ vanishes on the boundary of $Q_y$ since $\partial I^d\subset\Sigma_0$), and here the Poincaré constant is $k^{-2}C$, where $C$ depends only on the dimension $d$. Then, by putting all the inequalities together, restricting to $\Omega$ and integrating by parts as in Lemma \ref{poincare}, we get $||u_k||_{L^2(\Omega)}\leq C||f||_{L^2(\Omega)}$. Let us now consider an arbitrary weakly convergent subsequence (not relabeled) and its limit $w_{f,\Sigma_0,\Omega}$. It is easy to see that the pointwise value of this limit function depends only on the local behaviour of $f$. In fact, the key assumption $\partial I^d\subset\Sigma_0$ produces small cubes around each point $x\in\Omega$ which do not affect each other. So, if $f=\sum_i f_i I_{A_i}$ is piecewise constant (the pieces $A_i$ being disjoint open sets, for instance), it happens that for large $k$ the value of $u_k$ at $x\in A_i$ depends only on $f_i$. It turns out that, for a piecewise constant function $f$, it holds $w_{f,\Sigma_0,\Omega}=f w_{1,\Sigma_0,\Omega}$. It is indeed clear that in this case ($f=1$), since we are simply homogenizing the function $u_{1,\Sigma_0,I^d}$, the limit of the whole sequence $(u_k)_k$ exists, does not depend on the global geometry of $\Omega$, but it is a constant and it is the same constant as if there was $I^d$ instead of $\Omega$. Then the constant is easy to be computed and is the constant $c(\Sigma_0)$ appearing in the statement. It remains now just to show that the equality $w_{f,\Sigma_0,\Omega}=fc(\Sigma_0)$ is true for any $L^2$ function $f$. The convergence of the whole sequence will then follow easily by uniqueness of the limit of subsequences. To get the result for a generic $f$, just take a sequence $(f_n)_n$ of piecewise constant functions approaching it in $L^2$ and notice that
$$k^2u_{f,\Sigma^k,\Omega}=k^2u_{f_n,\Sigma^k,\Omega}+k^2u_{f-f_n,\Sigma^k,\Omega}.$$
The first term here weakly converges to $f_nc(\Sigma_0)$ as $k\tto+\infty$, while the second is bounded in the $L^2$ norm by $C||f-f_n||_{L^2(\Omega)}$. This means that any weak limit of subsequences of $(u_k)_k$ must be close in the $L^2$ norm to $f_nc(\Sigma_0)$, i.e.
$$||w_{f,\Sigma_0,\Omega}-f_nc(\Sigma_0)||_{L^2(\Omega)}\leq C||f-f_n||_{L^2(\Omega)},$$
which implies, letting $n\tto+\infty$, $w_{f,\Sigma_0,\Omega}=fc(\Sigma_0)$.
\end{proof}
Now we want to build efficient sets $\Sigma_0$ satisfying the key assumption of our previous Lemma, that is $\partial I^d\subset\Sigma_0$ (we will call those sets for which such an inclusion boundary-covering sets).

\begin{lem}\label{trovosigmabuono}
For any $\alpha>0$ and any $\ve>0$ there exists $n_0\in\N$ such that for any $n>n_0$ we find $\alpha''<\alpha$ and a set $\Sigma\in\A(\alpha'',n)(I^d)$ which is boundary-covering, with $n^{2/d}\int_{I^d}u_{1,\Sigma,I^d}d\lcal^d<(1+\ve)\theta^+(\alpha)$.
\end{lem}
\begin{proof}
Let us fix $\delta>0$ and $\alpha'<\alpha$ such that $\theta(\alpha')<(1+\delta)\theta^+(\alpha)$. By definition of $\theta(\alpha')$, we may find $\Sigma_1\in\A(\alpha',n_1)(I^d)$ such that 
$$n_1^{2/d}\int_{I^d}u_{1,\Sigma_1,I^d}\,d\lcal^d<(1+\delta)\theta(\alpha')$$
and, moreover, the number $n_1$ may be chosen as large as we want. Now we enlarge the set $\Sigma_1$ to get a new set $\Sigma_2$ which is boundary-covering: we add to $\Sigma_1$ some $m$ balls of radius $r=\alpha'n_1^{-1/d}$ (the same radius of the balls composing $\Sigma_1$). In order to cover $\partial I^d$ the number of balls we need does not exceed $C/r^{d-1}$, so we have $m\leq C(\alpha') n_1^{1-1/d}$. It is possible to choose $n_1$ so that 
$$m\leq \delta n_1\text{ and }\alpha'\left(\frac{n_1+m}{n_1}\right)^{1/d}=\alpha''<\alpha.$$
This is useful, since 
$$\Sigma_2\in\A\left(\alpha'\left(\frac{n_1+m}{n_1}\right)^{1/d}\vspace{-0.3cm},n_1+m\right)(I^d)=\A(\alpha'',n_2)(I^d),$$
where we set $n_2=n_1+m$. Moreover
$$n_2^{2/d}\int_{I^d}u_{1,\Sigma_2,I^d}\,d\lcal^d\leq \left(\frac{n_2}{n_1}\right)^{2/d}n_1^{2/d}\int_{I^d}u_{1,\Sigma_1,I^d}\,d\lcal^d<(1+\delta)^{2+2/d}\theta^+(\alpha).$$
Now, if we are given a large number $n$, we just need to homogenize the set $\Sigma_2$. By homogenization of order $k$ of a set $S\subset I^d$ into a domain $A$ we mean the set $A\cap\bigcup_{y\in\ k^{-1}\Z^d}y+k^{-1}S.$
Here we take the homogenization of order $k$ of $\Sigma_2$ into $I^d$, which is a set
$$\Sigma\in\A(\alpha'',k^dn_2)(I^d)\subset\A(\alpha'',n)(I^d),$$
where we choose $k$ such that $k^dn_2\leq n< (k+1)^d n_2$. For this set $\Sigma$ it holds
$(k^dn_2)^{2/d}\int_{I^d}u_{1,\Sigma,I^d}d\lcal^d=n_2^{2/d}\int_{I^d}u_{1,\Sigma_2,I^d}d\lcal^d,$ so that
$$n^{2/d}\int_{I^d}u_{1,\Sigma,I^d}d\lcal^d\leq \left(\frac{k+1}{k}\right)^{2/d}(1+\delta)^{2+2/d}\theta^+(\alpha).$$
If $n>n_2(1+\delta^{-1})^d$, then it holds $k>\delta^{-1}$ and $1+1/k<1+\delta$, so that we get
$$n^{2/d}\int_{I^d}u_{1,\Sigma,I^d}d\lcal^d\leq (1+\delta)^{2+4/d}\theta^+(\alpha).$$
It is now sufficient to choose $\delta$ sufficiently small so that $(1+\delta)^{2+4/d}<1+\ve$ and then set $n_0=n_2(1+\delta^{-1})^d$.
\end{proof}
We are now ready to start with the $\Gamma-\limsup$ main part.
We will start from a very particular class of measures. Let us call piecewise constant those probability measures $\mu\in\pical(\clom)$ which are of the form 
$$\mu=\rho\cdot\lcal^d,\text{ with }\rho\in L^1(\Omega),\,\int_{\Omega}\rho\, d\lcal^d=1,\,\rho>0,$$
for a piecewise constant function $\rho=\sum_{i=0}^m \rho_i I_{\Omega_i}$, the pieces $\Omega_i$ being disjoint Lipschitz open subsets with the possible exception of $\Omega_0=\Omega\setminus\bigcup_{i=1}^m\Omega_i$. To simplify the notation, let us also define the functional $\tilde{F}$, which is the same as $F$ with the only difference that we replace $\theta^-$ by $\theta^+$:
$$\tilde{F}(\mu)=\int_{\Omega}f^2\frac{\theta^+(\alpha\mu_a^{1/d})}{\mu_a^{2/d}}\,d\lcal^d=\int_{\Omega}f^2g_{\alpha}^+(\mu_a)\,d\lcal^d.$$

\begin{prop}
Under the same hypotheses of Theorem \ref{teoprinc}, it holds
$$ F^+(\mu)\leq \tilde{F}(\mu),\text{ where }F^+=\Gamma-\limsup_n F_n,$$
for any a piecewise constant measure $\mu\in\pical(\clom)$. This means that, for any such a measure $\mu$ and any $\ve>0$, there exists a sequence of sets $(\Sigma_n)_n$ such that $\mu_{\Sigma_n}$ weakly-* converges to $\mu$, $\Sigma_n\in\aan(\Omega)$ and moreover it holds 
$$\limsup_n n^{2/d}\int_{\Omega}fu_{f,\Sigma_n,\Omega}\,d\lcal^d\leq (1+\ve)\int_{\Omega}f^2\frac{\theta^+(\alpha\rho^{1/d})}{\rho^{2/d}}\,d\lcal^d.$$
\end{prop}
\begin{proof}
First of all, let us consider the numbers $\alpha\rho_i^{1/d}$, which will appear as arguments of the function $\theta^+$. We know by applying Lemma \ref{trovosigmabuono} to all of them that there exists a common number $n_0$ and some sets $\Sigma^i\in\A(\alpha''\rho_i^{1/d},n_0)(I^d)$, which are all boundary-covering and such that 
$$n_0^{2/d}\int_{I^d}u_{1,\Sigma^i,I^d}d\lcal^d<(1+\ve)\theta^+(\alpha\rho_i).$$
Now, if we are given some numbers $k(n,i)$, we define the sets $\Sigma_n^i$ by homogenizing into $\Omega_i$ the set $\Sigma_i$ of order $k(n,i)$  i.e.
$$\Sigma_n^i=\overline{\Omega_i}\cap\bigcup_{y\in A_{n,i}}y+k(n,i)^{-1}\Sigma^i,$$
where $A_{n,i}=\left\{y\in k(n,i)^{-1}\Z^d|(y+k(n,i)^{-1}\Sigma^i)\cap\Omega_i\neq\emptyset\right\}$. We define as well
$$\hat{\Sigma}_n^i=\bigcup_{y\in A_{n,i}}y+k(n,i)^{-1}\Sigma^i,$$
i.e. without intersecting the balls with $\overline{\Omega_i}$.

Then we choose $\Sigma_n=\bigcup_i\hat{\Sigma}_n^i\cup\tilde{\Sigma}_n$ where $\tilde{\Sigma}_n$ is a union of balls of radius $\alpha''n^{-1/d}$ covering the union of the boundaries $\partial\Omega_i$ inside the interior of $\Omega$. 
The number of balls we use in $\Sigma_n$ is approximatively $|\Omega_i|k(n,i)^d n_0$ in each zone $\Omega_i$, plus $Cn^{1-1/d}$ for the set $\tilde{\Sigma}_n$, where $C$ depends on $\alpha''$ and the total perimeter of the partition $(\Omega_i)_i$. In each zone $\Omega_i$ the radius of the balls is given by $\alpha''\rho_i^{1/d} n_0^{-1/d}/k(n,i)$. This means that the sequence $\Sigma_n$ we are building is admissible (i.e. $\Sigma_n\in\aan(\Omega)$ and $\mu_{\Sigma_n}\deb\mu$) if we have
\begin{gather*}
\frac{\alpha''\rho_i^{1/d}}{ n_0^{1/d}k(n,i)}\leq \alpha n^{-1/d} \text{ for $n$ large enough and for }i=0,\dots,m;\\
\sum_{i=0}^m |\Omega_i|k(n,i)^d n_0 + Cn^{1-1/d} \leq n \text{ and is asymptotic to }n;\\
\frac{k(n,i)^d n_0 }{n}\tto \rho_i \text{ for }i=0,\dots,m.
\end{gather*}
The first conditions involves radii, the second one the number of balls and the third one is related to weak convergence. All these conditions are satisfied if we set 
$$k(n,i)=\left\lfloor \left((1-Cn^{-1/d})\rho_i\frac{n}{n_0}\right)^{1/d}\right\rfloor.$$ 

We want now to estimate the values of the functionals on this sets $\Sigma_n$. We have used the set $\tilde{\Sigma}_n$ covering the internal boundaries of the sets $\Omega_i$ in order to get a local behaviour in which different zones $\Omega_i$ are independent on each other. The quantity we want to estimate is $$n^{2/d}\int_{\Omega}fu_{f,\Sigma_n,\Omega}\,d\lcal^d=\sum_{i=0}^m\int_{\Omega_i}fn^{2/d}u_{f,\Sigma_n,\Omega}\,d\lcal^d\leq \sum_{i=0}^m\int_{\Omega_i}fn^{2/d}u_{f,\Sigma_n^i,\Omega_i}\,d\lcal^d.$$
In the last inequality we have used the fact that $\tilde{\Sigma}_n$ allows adding Dirichlet boundary conditions on each $\Omega_i$ and $\hat{\Sigma}^i_n\supset\Sigma^i_n$. The disintegration of the integral here performed allows us applying Lemma \ref{convergenza a cf} on each $\Omega_i$. Notice first that $n\approx k(n,i)^d\rho_i^{-1}n_0$. Then, by Lemma \ref{convergenza a cf}, we know that $k(n,i)^2u_{f,\Sigma_n^i,\Omega_i}\deb c(\Sigma^i)f$, where the convergence is weak in $L^2$. By our choice of $\Sigma^i$ we have $c(\Sigma^i)<(1+\ve)n_0^{-2/d}\theta^+(\alpha\rho_i^{1/d})$, so we get
$$\lim_n n^{2/d}\int_{\Omega_i}fu_{f,\Sigma^i_n,\Omega_i}\,d\lcal^d<(1+\ve)\rho_i^{-2/d}\theta^+(\alpha\rho_i^{1/d})\int_{\Omega_i}f^2\,d\lcal^d,$$
and, summing up,
$$\lim_n n^{2/d}\int_{\Omega}fu_{f,\Sigma_n,\Omega}\,d\lcal^d\leq \int_{\Omega}f^2\frac{\theta^+(\alpha\rho^{1/d})}{\rho^{2/d}}\,d\lcal^d,$$
that is our thesis.
\end{proof}
We must extend our result to non piecewise constant measures and replace $\theta^+$ by $\theta^-$. To simplify the notation, we will use the function $g_{\alpha}$ defined in Section \ref{sec2}.
\begin{prop}\label{approx}
For any $\mu\in\pical(\clom)$ it holds 
$$F^+(\mu)\leq F(\mu)=\int_{\Omega}f^2 g_{\alpha}^-(\mu_a)\,d\lcal^d.$$
\end{prop}
\begin{proof}
Given $\mu=\rho\cdot\lcal^d$, with $\rho>c>0$, we may choose a sequence $\mu_n=\rho_n\cdot\lcal^d$ where $\rho_n\tto\rho$ a.e. and $\mu_n$ are piecewise constant with $\rho_n>c$. This may be done by approximating $\mu$ in $L^1$ first by regular functions, then by functions which are constant on cubes, for instance. So, by lower semicontinuity of $F^+$ (see \cite{introgammaconve}), we get
$$F^+(\mu)\leq \liminf_n \tilde{F}(\mu_n)\leq\int_{\Omega}f^2g_{\alpha}^+(\rho)\,d\lcal^d,$$
where we used the inequality $g_{\alpha}^+(\rho_n)\leq g_{\alpha}^+(c)$ to have dominated convergence of $g_{\alpha}^+(\rho_n)$ to their a.e. limit, which is estimated by upper semicontinuity by $g_{\alpha}^+(\rho)$. So we have extended our inequality to any absolutely continuous measure with positive density bounded away from $0$. 
Now take $\mu=\rho\cdot\lcal^d$ without the assumption $\rho>c$ and take $\rho_c=(\rho + c)1_{\rho\leq M}+(\rho-\ve_c)1_{\rho> M}$, where $\ve_c$ is chosen so that $\int_{\Omega}\rho_c\,d\lcal^d=1$. It is clear that $\ve_c\tto 0$ when $c\tto 0$ and that $\mu_c=\rho_c\cdot\lcal^d$ converges to $\mu$. So we have
\begin{eqnarray*}
F^+(\mu)&\leq& \liminf_{c\tto 0} \tilde{F}(\mu_n)=
\lim_{c\tto 0}\int_{\rho\leq M}f^2g_{\alpha}^+(\rho_c)\,d\lcal^d+\lim_{c\tto 0}\int_{\rho> M}f^2g_{\alpha}^+(\rho_c)\,d\lcal^d\\
&=&\int_{\rho\leq M}f^2g_{\alpha}^-(\rho)\,d\lcal^d+\int_{\rho> M}f^2g_{\alpha}^+(\rho)\,d\lcal^d\\
&\leq&\int_{\Omega}f^2g_{\alpha}^-(\rho)\,d\lcal^d+g_{\alpha}^+(M)\int_{\rho> M}f^2\,d\lcal^d.
\end{eqnarray*}
The number $M$ being arbitrary, we get easily $F^+(\mu)\leq \int_{\Omega}f^2g_{\alpha}^-(\rho)\,d\lcal^d$ by letting $M$ tend to $+\infty$. We need now only to get the result for measures of the form $\mu=\rho\cdot\lcal^d+\nu$ with a singular part $\nu$. Such a measure may be approximated by a sequence of measures $\mu_n=\rho 1_{A_n}+v1_{A_n^c}$ where the $v-$part carries the mass of the singular part $\nu$ and $1_{A_n}\tto 1_{\Omega}$ in $L^1$. We may choose $v$ so that it is $v>M$ with arbitrary large $M$, and in particular we may have $g_{\alpha}^-(v)=0$ since the function $\theta$ vanishes for large values of the argument. So we get
$$F^+(\mu)\leq\liminf_n F(\mu_n)=\liminf_n\int_{A_n}f^2g_{\alpha}^-(\rho)\,d\lcal^d\leq\int_{\Omega}f^2g_{\alpha}^-(\rho)\,d\lcal^d.$$
With this the proof is over.
\end{proof}

\section{Some properties of the function $\theta$}\label{sec4}

In this section we study some properties of the function $\theta$. We already know that it is decreasing and that it vanishes from a certain point on. Our first goal is to show the equality $\theta=\theta^-$, so that we can get rid of the l.s.c. envelopes.
\begin{prop}
For any $\alpha>0$ it holds $\theta(\alpha)=\theta^-(\alpha)$.
\end{prop}
\begin{proof}
From its definition:
\begin{eqnarray*}
\theta(\alpha)&=&\inf\left\{\liminf_n n^{2/d}F(\Sigma_n,1,Q)\left|\Sigma_n\in\aan(Q)\right.\right\}\\
&=&\liminf_n n^{2/d} \min\left\{F(\Sigma,1,I^d)\left|\Sigma\in\aan(I^d)\right.\right\}.
\end{eqnarray*}
It is then clear, by $\Gamma-$convergence, that such a $\liminf$ is in fact a limit equal to the minimum of the limit problem in $I^d$ with $f=1$. So 
$$\theta(\alpha)=\min\left\{ \int_{I^d}g_{\alpha}^-(\rho)\,d\lcal^d \left|\rho\geq 0,\,\rho\in L^1(I^d), \int_{I^d}\rho \,d\lcal^d=1 \right.\right\},$$
since the minimum will be certainly achieved on an absolutely continuous measure (as already mentioned in Theorem \ref{opticond}, as a consequence of $g_{\alpha}^-$ being decreasing). Moreover, it is clear that $\rho=1$ achieves the minimum: in fact, taken an optimal $\rho$, we may homogenize it so that we get a sequence $\rho_n$ weakly converging to the constant $1$. Each of these measures will be a minimizer since the functional only takes into account the measure of the level sets of $\rho_n$, which are the same as in $\rho$. So, by lower semicontinuity, $\rho=1$ will be a minimizer. We deduce
$\theta(\alpha)=g_{\alpha}^-(1)=\theta^-(\alpha),$ and this yields that $\theta$ is lower semicontinuous.
\end{proof}
Now we have $g_{\alpha}=g_{\alpha}^-$: our next step will be the following:
\begin{prop}
For any $\alpha>0$, the function $g_{\alpha}^-$ is convex.
\end{prop}
\begin{proof}
This is a consequence of the lower semicontinuity of $F$, which is implied by the fact that it is a $\Gamma-$limit (see \cite{introgammaconve}). This is quite standard (by necessary conditions on functionals of this kind in order to be l.s.c., see for instance \cite{bb2}), but not immediate since we are restricted to probability measures. Let us prove it, for the sake of completeness. In fact, given two values $x,\,y>0$ and $t\in[0,1]$ it is possible to build a probability measure $\mu$ on a domain $R=(0,a)\times I^{d-1}$ with density $x$ in $(0,at)\times I^{d-1}$ and $y$ in $(at,a)\times I^{d-1}$ ($a$ has to be chosen so that $a(tx+(1-t)y)=1$). By homogenizing $\mu$ on the first coordinate we get a sequence $\mu_n$ in which the value of the functional $F$ is constant and equal to $F(\mu)$ (we take $f=1$). Then we notice $\mu_n\deb\mu_{\infty}$, where $\mu_{\infty}$ is the measure with density $tx+(1-t)y$. By lower semicontinuity we get
$$ag^-(tx+(1-t)y)\leq atg^-(x)+a(1-t)g^-(y),$$
which gives the required convexity.
\end{proof}
\begin{cor}
The function $\theta$ is locally Lipschitz continuous on $(0,+\infty)$.
\end{cor}
\begin{proof}
This is a consequence of the equality $\theta(x)=g_1(x^d)x^2$, where both the factors in the right hand side product are locally Lipschitz functions of the variable $x$.
\end{proof}
The next thing to do with $\theta$ is to prove some estimate. In particular it is necessary to get estimates from below, so that $\theta$ is not identically $0$, otherwise our limit functional $F$ would be trivial.
%\begin{prop}\label{lower on theta}
%There exists a number $\alpha_0>0$ such that $\theta(\alpha_0)>0$.
%\end{prop}
%
%\begin{proof}
%Given $n\in\N$ and $\alpha_0>0$, let us divide $I^d$ into $k_n^d$ small cubes of side $1/k_n$, where we choose $k_n=\alpha_0^{-1}\left\lfloor n^{1/d} \right\rfloor$. When we take $\Sigma\in\A(\alpha_0,n)(I^d)$, the centres of the $n$ balls of $\Sigma$ will be placed in $n$ of such cubes, and, since the radii satisfy $r=\alpha_0 n^{-1/d}\leq 1/k_n$, the number of cubes intersected by the balls cannot exceed $3^dn$. So $u_{1,\Sigma,I^d}$ is greater than the function which vanishes on all these cubes and coincides with $v_i$ given by the solution of 
%$$\begin{cases}-\Delta v_i=1&\text{ in }Q_i\\
%							v_i=0&\text{ in }\partial Q_i,\end{cases}$$
%in any other cube $Q_i$ of this partition. It holds, by a scaling argument, $\int_{Q_i}v_i\,d\lcal^d=k_n^{-(2+d)}C(d)$.
%Consequently it holds 
%$$n^{2/d}\int_{I^d}u_{1,\Sigma,I^d}\,d\lcal^d\geq n^{2/d}\left(k_n^{2/d}-3^dn\right)k_n^{-(2+d)}C(d),$$
%and, passing to $\Sigma$ minimizing and then to the $\liminf$,
%$$\theta(\alpha_0)\geq C(d)\left(\frac{1}{\alpha_0^d}-3^d\right)\alpha_0^{d+2}.$$
%From this it is clear that for small $\alpha_0$ we will have $\theta(\alpha_0)>0$.
%\end{proof}
Anyway, we start from an estimate of $\theta$ from above:
\begin{prop}
There exists a constant $C=C(d)$ such that, for any $\alpha>0$, it holds $\theta(\alpha)\leq C\alpha^{2-d},$ for $d\geq 3$, or $\theta(\alpha)\leq C\log((\sqrt{2}\alpha)^{-1}),$ for $d=2$.
\end{prop}
\begin{proof}
It is sufficient to consider a particular sequence (or even a subsequence) of sets $\Sigma_n\in\aan(I^d)$ and then to compute the $\liminf$ in the definition of $\theta(\alpha)$. We choose to consider just the numbers $n$ of the form $n=k^d$, and to build, for each $k\in\N$, a set $\Sigma_n$ which is composed by $n=k^d$ balls of radius $\alpha/k$, with their centres placed at the middle points of the $k^d$ cubes of side $1/k$ of a regular lattice partitioning the cube $I^d$. First, we notice that it holds $u_{1,\Sigma_n,I^d}\leq v_n$, where $v_n$ is the solution to the problem
$$\begin{cases}-\Delta v_n=1&\text{ in }I^d\setminus\Sigma_n\\
							v_n=0&\text{ in }\Sigma_n,\\
							\frac{\partial}{\partial n}v_n=0&\text{ on }\partial I^d.\end{cases}$$
The inequality between the solutions of Dirichlet and Neumann problems comes from maximum principle.
By scaling arguments, it is clear that the energy $n^{2/d}\int_{I^d}v_n\,d\lcal^d$ equals the energy $\int_{I^d}v_1\,d\lcal^d$. If we set $r_0=\sqrt{d}/2$, which is the radius of the smallest ball containing the cube $I^d$ and centred at its same centre, it holds $v_1\leq w$, where $w$ is the solution of the Neumann problem on such a ball:
$$\begin{cases}-\Delta w=1&\text{ in }B(x_0,r_0)\setminus \ball{x_0}{\alpha}\\
							w=0&\text{ in }\ball{x_0}{\alpha},\\
							\frac{\partial}{\partial n}w=0&\text{ on }\partial B(x_0,r_0),\end{cases}$$
where $x_0$ is the centre of the cube $I^d$.
This solution may be explicitly computed, it is radially symmetric and is given by
\begin{eqnarray*}
w(x)&=&k(d)(\alpha^{2-d}-r^{2-d})-\frac{1}{2d}(r^2-\alpha^2)\text{  with }r=|x-x_0|,\text{ if }d\geq 3;\\
w(x)&=&k\log\left(\frac{r}{\alpha}\right)-\frac{1}{4}(r^2-\alpha^2)\text{  with }r=|x-x_0|,\text{ if }d=2.
\end{eqnarray*}
It turns out that 
$$\frac{\partial}{\partial n}w>0\text{ on }\partial I^d,$$
which allows using maximum principle to get $v_1\leq w$. So it is sufficient to compute the integral of $w$ on $I^d$, which can be estimated by $k(d)\alpha^{2-d}$ in the case $d\geq 3$, and by $k\log(\alpha^{-1})+k\log r_0= k\log((\sqrt{2}\alpha)^{-1})$ in the case $d=2$.
\end{proof}

A similar estimate from the other side may be obtained, as we show in our next proposition. Here the techniques we use are quite different and much more related to shape optimization and PDEs.

\begin{prop}
For $\alpha<t_1$ the following estimates hold:
\begin{eqnarray}
-\theta'(\alpha)&\geq& \frac{\alpha^{1-d}}{d\omega_d}-2\frac{\alpha}{d}\\
\theta(\alpha)&\geq& \frac{\alpha^{2-d}}{d(d-2)\omega_d}-C\text{ if }d\geq 3,\\
\theta(\alpha)&\geq& \frac{\log(\alpha^{-1})}{2\pi}-C\text{ if }d=2,
\end{eqnarray}
for a suitable constant $C=C(d)$, where $\omega_d$ denotes the measure of the $\R^d-$unit ball.
\end{prop}
\begin{proof}
For a fixed $n\in\N$ and fixed points $(x_i)_{i=1,\dots,n}\in I^d$, let us consider the sets $\Sigma_{\alpha}=\bigcup_{i=1}^n\ball{x_i}{\alpha n^{-1/d}}\in\aan(I^d)$ and $\Omega_{\alpha}=I^d\setminus\Sigma_{\alpha}$. The following estimate holds, by Holder inequality:
\begin{equation}\label{byholder}
\haus^{d-1}(\partial\Omega_{\alpha})\left(\int_{\partial\Omega_{\alpha}}|\frac{\partial}{\partial n}u_{1,\Sigma_{\alpha},I^d}|^2\,d\haus^{d-1}\right)\geq\left(\int_{\partial\Omega_{\alpha}}\frac{\partial}{\partial n}u_{1,\Sigma_{\alpha},I^d}\,d\haus^{d-1}\right)^2=|\Omega_{\alpha}|^2,
\end{equation}
where the last equality comes from integrating by parts $\int_{\Omega_{\alpha}}-\Delta u_{1,\Sigma_{\alpha},I^d}\,d\lcal^d$. Notice that it holds as well
\begin{equation}\label{hausd-1}
\haus^{d-1}(\partial\Omega_{\alpha})\leq nd\omega_d(\alpha n^{-1/d})^{d-1}=d\omega_d\alpha^{d-1}n^{1/d}.
\end{equation}
Then, we recall that, by using shape derivative (\cite{HP}), it holds 
\begin{equation}\label{shapeder}
-\frac{d}{d\alpha} F(\Sigma_{\alpha},1,I^d)=n^{-1/d}\int_{\partial\Omega_{\alpha}}|\frac{\partial}{\partial n}u_{1,\Sigma_{\alpha},I^d}|^2\,d\haus^{d-1},
\end{equation}
since we perturbate $\Omega_{\alpha}$ by a vector field which is normal to the boundary of the balls and proportional to $n^{-1/d}$. So we have, by putting together \eqref{byholder}, \eqref{hausd-1} and \eqref{shapeder}
$$-\frac{d}{d\alpha} n^{2/d}F(\Sigma_{\alpha},1,I^d)\geq \frac{\alpha^{d-1}|\Omega_{\alpha}|^2}{d\omega_d}\geq\frac{\alpha^{d-1}(1-2\omega_d\alpha^d)}{d\omega_d},
$$
where in the last inequality we have used $|\Omega_{\alpha}|^2\geq (1-\omega_d\alpha^d)^2\geq 1-2\omega_d\alpha^d $.
So, by integrating over an arbitrary interval $(\alpha_1,\alpha_0)\subset(0,t_1)$, we have
\begin{equation}\label{integrando}
n^{2/d}F(\Sigma_{\alpha_1},1,I^d)\geq n^{2/d}F(\Sigma_{\alpha_0},1,I^d)+\frac{1}{d\omega_d}\int_{\alpha_1}^{\alpha_0}\alpha^{1-d}(1-2\omega_d\alpha^d)\,d\alpha.
\end{equation}
We now compute explicitly the right hand side of \eqref{integrando} in the case $d\geq 3$, getting
$$n^{2/d}F(\Sigma_{\alpha_1},1,I^d)\geq
n^{2/d}F(\Sigma_{\alpha_0},1,I^d)+\frac{\alpha_1^{2-d}-\alpha_0^{2-d}}{d(d-2)\omega_d}-\frac{\alpha_0^2-\alpha_1^2}{d}.$$
Then, passing to the $\inf$ over $(x_i)_i$ and to the $\liminf$ over $n$, we get
\begin{equation}\label{cond3opiu}
\theta(\alpha_1)\geq \theta (\alpha_0)+\frac{\alpha_1^{2-d}-\alpha_0^{2-d}}{d(d-2)\omega_d}-\frac{\alpha_0^2-\alpha_1^2}{d}.
\end{equation}
This gives the estimate on $\theta$ we were looking for, with 
$$C(d)=t_1^2/d+t_1^{2-d}/[d(d-2)\omega_d].$$
Anyway, from \eqref{cond3opiu}, we can also infer the estimate on $\theta'$: we let $\alpha_1\tto\alpha_0^-$ and divide by $\alpha_0-\alpha_1$. On those points $\alpha_0$ where $\theta$ is differentiable (i.e. almost everywhere, since $\theta$ is locally Lipschitz) it holds exactly
$$-\theta'(\alpha_0)\geq \frac{\alpha_0^{1-d}}{d\omega_d}-2\frac{\alpha_0}{d}.$$
In the case $d=2$ it is sufficient to compute again the integral in \eqref{integrando}, getting
$$nF(\Sigma_{\alpha_1},1,I^d)\geq nF(\Sigma_{\alpha_0},1,I^d)+\frac{\log(\alpha_1^{-1})-\log(\alpha_0^{-1})}{2\pi}-\frac{\alpha_0^2-\alpha_1^2}{2}.$$
Then the conclusion follows in the same way. Here $C=\log(t_1^{-1})/2\pi + t_1^2/2.$
\end{proof}

We summarize now all the results on $\theta$ we have got in this section.
\begin{teo}
The function $\theta:]0,+\infty[\tto [0,+\infty[$ is a nonincreasing function, locally Lipschitz on $]0,+\infty[$, such that $x\mapsto \theta(x)x^{-2/d}$ is a convex function. Moreover $\theta$ is not identically $0$ but it vanishes from a certain point on, i.e. $\theta(x)=0$ for any $x\geq t_1$ with $t_1\leq \sqrt{d}/2$. Finally it holds 
\begin{gather*}
C_1x^{2-d}-C\leq\theta(x)\leq C_2x^{2-d}\mbox{ if }d\geq 3\\
C_1|\log x|-C\leq\theta(x)\leq C_2|\log x|\mbox{ if }d=2.
\end{gather*}
\end{teo}

We terminate this section by stressing the interest in finding explicit minimizing sequences for the case $f=1$ and $\Omega=I^d$, since this would give the value of $\theta$. In analogy to what happens in the location problem (see \cite{hexag} and \cite{toth}), we may conjecture that minimizing sequences are given by placing the centres of the balls on some kind of regular grids and, moreover, as far as $d=2$, hexagonal regions with balls in the middle seem to be good candidates. This would lead to a better knowledge of $\theta$, thus letting us get better information on the minimizer $\mu$ for the limit functional $F$ by Theorem \ref{opticond}. As a weaker conjecture, we may think that the function $\theta$, for which we have proven upper and lower estimates in term of $x^{2-d}$ or $|\log(x)|$, is such that there also exist the limits 
$$\lim_{x\tto 0^+}\frac{\theta(x)}{x^{2-d}}\text{ for }d\geq 3\text{ or }\lim_{x\tto 0^+}\frac{\theta(x)}{|\log(x)|}\text{ for }d=2.$$
In fact what we have already proven is just that the ratios above are bounded. For instance, in the case $d\geq 3$, should $\theta(x)$ actually behave like $x^{2-d}$ for $x$ near $0$, we could derive, for small $\alpha$, a behaviour like
$$\mu_{\alpha}\approx cf.$$
%even if we need some results also on the derivatives of $\theta$ and $g_{\alpha}$ to get a similar result.

\section{The one dimensional case}\label{sec5} 

In the case of dimension $1$ we are able to compute the function $\theta$ explicitly.

Everything is, in fact, simpler in dimension $1$, since the balls we remove are intervals which disconnect the domain of the differential equation (an ODE in this case), and so we can compute explicitly the solution. We have already pointed out that in dimension $1$ the compliance optimization problem is well-posed also for finite unions of points (i.e. the case $\alpha=0$), and not only for small intervals. So we will take into account also the value of $\theta(0)$.
\begin{teo}\label{theta in dim1}
For any $\alpha\geq 0$, if $d=1$, it holds
$$\theta(\alpha)=\begin{cases}\frac{1}{12}(1-2\alpha)^3&\text{ if }\alpha\leq \frac{1}{2}\\
                              0                       &\text{ if }\alpha\geq \frac{1}{2}.\end{cases}$$
\end{teo}
\begin{proof}
Given an interval $J$ whose length is $l>0$, the solution of the problem 
$$\begin{cases}-u''=1&\text{ in }J\\
                  u=0&\text{ on }\partial J\end{cases}$$
is given by $u(x)=x(l-x)/2$ (if we suppose $J=(0,l)$), and so its integral on $J$ is $l^3/12$. When we are given $\alpha\in[0,1]$, we take into consideration disjoint union of $n$ intervals, for a total length of $2\alpha$ (in this case it is necessary to have disjoint intervals if we want to minimize compliance). So the energy of the configuration is
$$\sum_{i=1}^n\frac{l_i^3}{12}\text{ under the condition }\sum_{i=1}^n l_i=1-2\alpha.$$
By convexity of $l\mapsto l^3$, the minimum of such a quantity is achieved by $l_i=n^{-1}(1-2\alpha)$, and so it holds
$$\theta(\alpha)=\liminf_n n^2 n \frac{(1-2\alpha)^3}{12n^3}=\frac{1}{12}(1-2\alpha)^3.$$
The case $\alpha>1/2$ follows trivially from $\theta(1/2)=0$ and $\theta$ being non negative and decreasing.
\end{proof}

As a consequence, the limit problem in dimension $1$ is completely known. Up to inverting the resulting function $g_{\alpha}'$ it is possible to find explicitly the minimizer, thanks to the expression for $\mu$ given by Theorem \ref{opticond}. 

We can also consider the sequence of minimization problems with $n$ points instead of $n$ balls: we have to consider
$$F_n(\mu)=\begin{cases}n^2 F(\Sigma,f,I)&\text{ if }\mu=\mu_{\Sigma}\text{ and }\sharp\Sigma\leq n\\
                        +\infty          &\text{ otherwise. }\end{cases}$$

\begin{teo}
Given $J=[a,b]$ and $f\in L^2(J)$, the sequence of functionals $(F_n)_n$ over $\pical(J)$ $\Gamma-$converges, with respect to weak* convergence on $\pical(J)$, to the functional $F$ given by
$$F(\mu)=\frac{1}{12}\int_J\frac{f^2}{\mu_a^2}\,d\lcal^1.$$
Moreover, $F$ has a unique minimizer $\mu_{opt}$, given by 
$$\mu_{opt}=c f^{2/3}\cdot\lcal^1, \text{ where }c=\left(\int_J f^{2/3}\,d\lcal^1\right)^{-1}.$$
\end{teo}
\begin{proof}
This statement follows by slight modifications of what proven in last section. The only point where $\alpha>0$ was actually used was Lemma \ref{poincare}. Anyway, in the one dimensional case, the inequality holds if we replace $|\left\{v=0\right\}|\geq \ve_0|I^d|$ by $\left\{v=0\right\}\neq\emptyset$. This allows performing again the proof, by following the same steps. To find explicitly the minimizer it is sufficient to use Lagrange multipliers, getting that $f^2/\mu_a^3$ is constant.
\end{proof}

\bigskip{\bf Acknowledgements.} The work on this subject started while the third author held a post-doc position at University of Pisa. Such a a post-doc position has been financially supported by the European Research Training Network {\it``Homogenization and Multiple Scales''(HMS2000)}. Grateful acknowledgements go to the Network and to the Department of Mathematics of Pisa for the hospitality.

\end{document}